\documentclass[11pt]{amsart}
\usepackage{latexsym}
\usepackage{mathrsfs,amsmath,amssymb,amscd,amsthm,amsfonts,verbatim,amsopn,color,graphics,eufrak}
\usepackage[all]{xy}
\usepackage[dvips]{epsfig} 
\def\PL{{\mathrm{\textup{Plim}}\,}}
\def\HC{{\mathrm{\textup{hocolim}}\,}}
\def\r{{\hat{\rho}}}
\def\d{{\mathcal{D}_{id}}}
\def\dd{{\mathcal{D}_\rho}}
\def\W{{W_\rho}}
\def\V{{\mathscr{V}_\rho}}
\def\L{{\Lambda_\rho}}

\def\C{{\mathbb C}}

\def\R{{\mathbb R}}

\def\PP{{\mathcal P}}
\def\AA{{\mathcal A}}

\def\CC{{\mathcal C}}
\def\DD{{\mathcal D}}
\def\FF{{\mathcal F}}

\def\II{{\mathcal I}}

\def\MM{{\mathcal M}}

\def\QQ{{\mathcal Q}}

\def\Ar{{\AA_\R}}
\def\Ac{{\AA}}
\def\S{{ Sal(\Ac)}}
\def\M{{\MM(\Ac)}}

\newtheorem{thm}{Theorem}[section]
\newtheorem{df} [thm]{Definition}
\newtheorem{lm}  [thm]{Lemma}
\newtheorem{crl} [thm]{Corollary}
\newtheorem{prop}[thm]{Proposition}

\newtheorem{rem}[thm]{Remark}

\newtheorem{expl}[thm]{Example}
\newtheorem{cstr}[thm]{Construction}
\numberwithin{equation}{section}

\newenvironment{pf}{\noindent {\bf Proof. }}{\hfill $\Box$\vspace{0.3cm}}

\begin{document}

\title[Diagram Models]
{Diagram models for the covers \\ of the Salvetti complex}

\author{Emanuele Delucchi}

\thanks{The author acknowledges support for this project by
        ETH research grant TH-10/02-3.}

\address{
Department of Mathematics, ETH Zurich, 8092 Zurich, Switzerland}
\email{delucchi@math.ethz.ch}

\begin{abstract}
  To every affine real arrangement of hyperplanes $\Ar$ we
  associate a family of diagrams of spaces over the face poset of the arrangement. We
  show that any cover of the complement of the complexification of
  $\Ar$ is homotopy equivalent to the homotopy colimit of one of
  the diagrams. More precisely, we show that any cover of the Salvetti
  complex is isomorphic to the order complex of the poset limit of one
  of the diagrams. We thus obtain explicit simplicial models for
  covers of the Salvetti complex. 
%A simplicial model of each cover is
%  explicitely presented as an order complex of a poset.
\end{abstract}

%%%%%%%%%%%%%%%%%%%%%%%%%%%%%%%%%%%%%%%%%%%%%%%%%%%%%%%%%%%%%%%%%%%%%%%%%
%%%%%%%%%%%%%%%%%%%%%%%%%%%%%%%%%%%%%%%%%%%%%%%%%%%%%%%%%%%%%%%%%%%%%%%%%
%%%%%%%%%%%%%%%%%%%%%%%%%%%%%%%%%%%%%%%%%%%%%%%%%%%%%%%%%%%%%%%%%%%%%%%%%

\maketitle

%%%%%%%%%%%%%%%%%%%%%%%%%%%%%%%%%%%%%%%%%%%%%%%%%%%%%%%%%%%%%%%%%%%%%%%%%
%%%%%%%%%%%%%%%%%%%%% Introduction     %%%%%%%%%%%%%%%%%%%%%%%%%%%%%%%%%%
%%%%%%%%%%%%%%%%%%%%%%%%%%%%%%%%%%%%%%%%%%%%%%%%%%%%%%%%%%%%%%%%%%%%%%%%%

\section*{Introduction}
\label{sect_intr}

%The aim of combinatorial models for the complement of an 
% arrangement is to clarify to what extent the
%combinatorics of the arrangement determines the topology of the
%complement.
%A famous problem is whether asphericity of the arrangement
%complement is a combinatorial property (the so-called $K(\pi,1)$
%problem for arrangements). This is still open even for the case of a
% complexified arrangement (i.e., an arrangement of complex hyperplanes
% with defining equations having only real coefficients), that is the
% type of arrangement which we model with our construction. 
%Furthermore, homological computations on covers are of interest (see
%for instance \cite{C}, where local system cohomology of some covers is related
%to the Milnor Fiber of the arrangement), and our model fits this
% setting too, since each cover of the complement can be modeled a a
% simplicial complex.\\

Let $V$ be a $d$-dimensional complex vector space. An {\em arrangement
  of hyperplanes} in $V$ is a finite set
$\mathcal{A}=\{H_i\}_{i=1,\dots ,n}$ of affine or linear codimension 1
subspaces of $V$. %From a combinatorial point of view, the union of
                  %the $H_i$ determines a stratification of $V$: the
                  %{\em combinatorial data} of the arrangement is the
                  %partial order of the strata, ordered by inclusion. 
The arrangement induces a stratification of the ambient space by its
hyperplanes and their intersections. The poset of strata, i.e. the
intersections in $\AA$, is customarily perceived as the
{\em combinatorial data} of the arrangement.
On the topological side, it is interesting to study the link $\bigcup H_i$ and the complement of the arrangement $\mathcal{M}(\mathcal{A}):=V\setminus\; \bigcup H_i$.

One of the main questions in arrangement theory  is to clarify to what
extent the combinatorial data of the arrangement determine topological
invariants of the complement or of the link of the arrangement.

Our interest restricts now to the case where $V$ is a complex vector
space. A  famous open question in this direction is the so-called $K(\pi ,1)$-problem.
An arrangement is said to be $K(\pi,1)$ if the higher homotopy groups of the complement $\pi_i(\M)$ vanish for $i>1$.
There are two large classes of arrangements that were shown to be
$K(\pi ,1)$ in classical works by Deligne (complexified simplicial
arrangements, see \cite{D}) and Falk-Randell and Terao (supersolvable arrangements, see
\cite{FR}, \cite{T}): both these classes admit a purely combinatorial characterization.
It is an open question whether being $K(\pi,1)$ is a combinatorial
property in general.

A very useful tool in studying the topology of arrangements are
combinatorial models for $\M$, i.e. cell complexes that are built from
the combinatorial data of the arrangement and that model $\M$ up to
homotopy equivalence or even homeomorphism. Many combinatorial models
for $\M$ exist at present. Among them let us only mention the {\em
  Salvetti complex}, introduced by Mario Salvetti in \cite{S}; it models the homotopy type of the complement of {\em complexified} arrangements (i.e., complex arrangements where the defining equations of the $H_i$ have real coefficients), and the model introduced by Bj\"orner and Ziegler in \cite{BZ}, that generalizes the work of Salvetti to any complex arrangement.

These combinatorial models led to important achievements in various
contexts, among others the computation of cohomology algebras of
complements of certain classes of subspace arrangements in purely
combinatorial terms, see \cite{FZ}, \cite{dL}, \cite{dLS}.
% such as the computation of the cohomology groups $H^{*} (\M)$ in purely combinatorial terms (in some extent even for the case of $subspace$ arrangements, i.e., when not only hyperplanes but any subspace can be element of $\AA$).\\

In this work we present a combinatorial model for the covers of the complement $\M$ of a complexified arrangement.

We hope that exploiting these models will bring some progress on the
$K(\pi,1)$-question (e.g. through explicit computations on the
universal cover), and will possibly simplify the computation of local system cohomology on $\M$.\\

In the first part of the article we will review the Salvetti complex,
oriented systems as introduced by Paris \cite{Pa1}, and diagrams of spaces. these three gadgets constitute the ingredients for the
formulation of our models.\\

In the second part we define certain diagrams of spaces, and we prove that these diagrams, in fact, serve as models for covers of arrangement complements.

In Section 4 we give the definition of the diagram belonging
to a given cover of an oriented system $(\Gamma,\sim)$ associated to a real arrangement $\AA$.
%The notion of oriented system was introduced by Paris, see
%\cite{Pa1}. 
We note that homotopy
colimits of our diagrams can always be written as an order complex of
an explicitly given poset and we prove that the homotopy
colimit of the diagram associated to the trivial cover of
$(\Gamma,\sim)$ is the Salvetti complex.
We then use this diagram to answer a question due to Michael Falk, who
asked for a generalization of a construction described in the proof
of lemma 2.3 of \cite{F}.

The Section 5 is entirely dedicated to proving that the
homotopy colimit of the diagram associated to a cover of
$(\Gamma,\sim)$ is a cover of the Salvetti complex in the topological
sense.

The following Section 6 introduces a formulation of the homotopy colimit
of our diagrams as CW complexes that is useful for considerations
concerning the fundamental group. We show that the fundamental group of
the homotopy colimit of a diagram is isomorphic to the fundamental group of the
oriented system to which the diagram is associated.

In the last Section 7 we summarize our work and we prove the main theorem
on diagram models: each cover $S'$ of the Salvetti complex is homotopy equivalent
to the homotopy colimit of a certain diagram model $\DD$. The order complex
of the poset limit of $\DD$ is isomorphic to $S'$ as a cover. As a
corollary we find that the homotopy colimit of the diagram
associated to the universal cover (in the oriented-system sense) of
$(\Gamma,\sim)$ gives the universal cover (in the topological sense)
of the Salvetti complex.

We conclude the article by illustrating the construction on a
simple but instructive example.

%%%%%%%%%%%%%%%%%%%%%%%%%%%INDICE%%%%%%%%%%%%%
\newpage
\tableofcontents

%%%%%%%%%%%%%%%%%%%%%%%%%%%%

\part{The toolkit.}
\label{toolkit}
In the first part of our article, we review the concepts and
definitions that will be used to formulate our model. We do
not aim for completeness: extensive references on each topic will be given. For a comprehensive introduction to arrangement theory we refer to \cite{OT}.

\section*{0. Basic Notions}

Let $\AA$ be an arrangement in $\R^d$,
i.e., a set of linear or affine codimension 1 subspaces of $\R^d$.

For an affine subspace $X$ of $\R^d$ it is standard to define a
subarrangement  $\mathcal{A}_X := \{H \in \mathcal{A} \;\vert\;
X\subseteq H\}$ of $\AA$, and the {\em restriction} of $\mathcal{A}$
to $X$, $\mathcal{A}^X := \{ X \cap H \;\vert\; H \in
\mathcal{A}\setminus\mathcal{A}_X\} $. The set of connected components
of $V\setminus\mathcal{A}$ we will denote with $\CC(\AA)$, the {\em chambers} of $\AA$.

There are two posets associated with $\mathcal{A}$ that we will
frequently use: the {\em intersection lattice} of $\AA$,
\begin{displaymath}
\mathcal{L}(\mathcal{A}):= \Big\{ X:= \bigcap_{\scriptscriptstyle H\in\mathcal{H}} H \;\Big\vert\; \mathcal{H} \subseteq \mathcal{A}, X \neq \emptyset \Big\},
\quad X\geq Y \iff X\subseteq Y,
\end{displaymath}
i.e., intersections of hyperplanes ordered by reversed inclusion,
and the {\em face poset} of $\AA$,
\begin{displaymath}
\mathcal{F}(\mathcal{A}):= \bigcup_{\scriptscriptstyle X\in\mathscr{L}(\mathcal{A})}\big\{ F \;\big\vert\; F \in \mathcal{C}(\mathcal{A}^{\scriptstyle X})\big\},
\quad  F_1 \geq F_2 \iff \overline{F_1}\subseteq\overline{F_2},
\end{displaymath}
i.e., cells in $\R^d$ that are chambers of restricted arrangements,
ordered by containment.

In the latter definition,  $\overline{F}$ denotes the closure of $F$
in the standard topology of $\R^d$. In fact, $\mathcal{F}(\AA)$ is the
poset of strata for the stratification of $\R^d$ into relatively open
polyhedral cells induced by the arrangement.
%Given a face $F$, $|F|$ will denote its affine hull.

Given a poset $\mathcal{P}$, it is standard to write $\Delta
(\mathcal{P})$ for the {\em order complex} of $\mathcal{P}$, i.e., the
abstract simplicial complex whose vertices are the elements of
$\mathcal{P}$ and whose simplices are the chains $p_0 > \dots > p_k$
in $\mathcal{P}$.

Throughout this paper, $\AA$ denotes a linear or affine {\em real} hyperplane arrangement. %$\AA_{\C}$ its complexification and $d$ the dimension of the ambient space.
Taking the same linear forms as for $\AA$ in complex $d$-space, we
obtain the complexification $\AA_{\C}$ of $\AA$.
We will denote by $\M$ the complement of $\AA_{\mathbb{C}}$ in $\C^d$, $\M:=\C^d\setminus\AA_{\C}$.

\section{The Salvetti Complex}
\label{salvetti}
A central tool in the study of the topology of $\M$ is the simplicial
complex $\S$, introduced by Salvetti in \cite{S}. We will recall its
definition and thereby essentially follow the original paper, although
we will adapt the notation for our purposes.

\subsection{Faces act on chambers}\label{fenincasina}
Before starting with the construction, we want to define a composition
of faces. Fix a face $F\in\FF(\AA)$ and let $\vert F\vert$ denote its affine hull.
The canonical projection $V\rightarrow V/ |F|$ induces a map
\begin{displaymath}
\pi_F : \mathcal{F}(\mathcal{A}) \mapsto \mathcal{F}(\mathcal{A}_{|F|}).
\end{displaymath}

For $P\in\mathcal{F}(\mathcal{A})$, let $F\circ P$ be the face of smallest dimension in $\pi_F^{-1}(\pi_F(P))$ whose closure contains $F$.

%\setlength{\unitlength}{3cm}
%\begin{picture}(2,2)
%\put(0,0){\line(5,3){2}}
%\put(0,0.8){\line(6,-1){2}}
%\put(0.7,0){\line(-1,5){0.3}}
%\put(0.2,1.5){\line(5,-2){1.7}}
%\end{picture}

It is not hard to see that this composition satisfies the following properties:
\begin{displaymath}\begin{array}{cl}
(\circ 1) &    F_1\circ F_2=F_2 \quad\textrm{for } \; F_1 \geq F_2
\;\textrm{in} \; \FF(\AA). \\
(\circ 2) & F\circ C \in \CC(\AA) \quad\textrm{ for any } F\in\FF(\AA), \; C\in\CC(\AA).\\
\end{array}\end{displaymath}
Writing $C_F$ for $F\circ C$ we can also state:
\begin{displaymath}\begin{array}{cl}
(\circ 3) &  (C_{F_1})_{F_2}=C_{F_1 \circ F_2}=C_{F_2} \quad\textrm{for } F_1 \geq F_2 \textrm{ and } C \in \CC(\AA). 
\end{array}\end{displaymath}
After these (perhaps at this point obscure) remarks, we move on to
define the Salvetti complex.

\subsection{Definition of $\S$}\label{defsal}

We start with the real arrangement $\AA$, and we choose a point $w(F)$ in the
interior of every face $F\in\FF(\AA)$. We write the complexification of $\R^n$ as $\C^n=\R^n\oplus\R^n$, and consider points
\begin{displaymath}
w(F,C):=w(F) \oplus [ w(F)- w(C_F)],
\end{displaymath}
for any $F\in\FF(\AA)$ and $C\in\CC(\AA)$.

\vspace{0.5cm}
Notice that, although indexed by $F$ and $C$, the points $w(F,C)$
actually depend on $F$ and $C_F$; only the relative position of $C$
with respect to $F$ is important.

These points will form the vertex set of the simplicial complex
$\S_{simpl}$. %The higher dimensional simplices will be convex hulls of
%some of the $w(F,C)$: 
The convex hull of $\{ w(F_0,C_0), \dots ,
w(F_d,C_d)\}$ is a simplex in $\S_{simpl}$ if and only if $F_0>F_1>\dots >F_d$
is a chain in $\FF$ and for any $i$, $(C_0)_{F_i}=(C_i)_{F_i}$. In
other words, each $(d+1)$-chain $\phi$ in $\FF$ and every chamber $C\leq max(\phi)$ determine a $d$-simplex $\Delta(\phi , C)$, and all simplexes are of this form.

There is another way to look at the Salvetti complex, which may be
more suggestive for some purposes. Indeed, one may take a vertex for
every chamber, and connect each pair of vertices corresponding to adjacent
chambers  by a pair of opposite oriented edges.  One attaches higher
dimensional cells to the resulting (oriented) graph in such a way that
 as many $d$-cells  correspond to every $(n-d)$-face $F$ of $\Ar$  as
there are chambers adjacent to $F$: The $d$-cell corresponding to the
face $F$ and the chamber $C$ is denoted by $[F,C]$. To see how it is
attached to the $(d-1)$-skeleton, consider the chamber $\tilde{C}$
opposite to $C$ with respect to $F$ and the positive directed minimal paths
from $C$ to $\tilde{C}$: these form the $1$-skeleton $[F,C]^1$ of
$[F,C]$, which consists of all cells of dimension $k\leq d$ whose
$1$-skeleton lies in $[F,C]^1$. This way one can construct a
CW-complex $\S_{CW}$ by induction on the skeleta. It is easy to see
that $\S_{simpl}$ is the barycentric  subdivision of $\S_{CW}$: in
fact, the simplices $\Delta(\phi,C)$ with $max(\phi)=F$ are simplices
of the baricentric subdivision contained in the closed cell $[F,C]$.  The two complexes are canonically homeomorphic, and we shall therefore identify them and use the name $\S$ for both constructions.

\subsection{The fundamental group $\pi_1(\S)$} 

Salvetti gave different presentations of $\pi_1(\S)$ in part 2 of
\cite{S}, and we refer there for explicit results. Here we only want
to remind that we can write the fundamental group as
\begin{equation}\label{pi}
\pi_1(\S) = \pi_1(\S^1) / \sigma ,
\end{equation}

where $\sigma$ denotes the smallest normal subgroup containing
$\{\partial [F,C] \vert \, codim(F)=2\}$, the set of boundaries of
$2$-cells (see e.g. \cite{M} ch. 7, Thm 2.1). In fact, the fundamental group is
determined by the $2$-skeleton of the complex.

From the construction of $\S$ as a CW-complex one can easily see
that $\sigma$ is generated by the elements $\gamma_1\gamma_2^{-1}$,
where the $\gamma_i$ are directed, minimal paths with same beginning
and end.

This last remark shows how  orientation of the edges carries the whole
information about the attaching of higher cells, and in particular
about the fundamental group. 
%At this point it is very natural to go to the next section, and introduce

%%%%%%%%%%%%%%%%%%%%%%%%% ORIENTED SYSTEMS %%%%%%%%%%%%%%%%%%%%%%%%%%%%%%%%%%%5

\section{Oriented Systems}

We turn our attention to a tool introduced by Paris \cite{Pa1}, which will
encode all the information of the arrangement that we need to define our diagram model.
\subsection{Definition of an oriented system $(\Gamma, \sim  )$}

Our starting point will be an oriented graph, i.e., a set of vertices $V$ and a set of edges $E$ consisting of some
ordered pairs $e=(v_1,v_2)$ of elements $v_i \in V$. A path $\alpha$
in an oriented graph is a sequence of edges
$e_1e_2\dots e_n$, each with an exponent $\epsilon_i =\pm 1$, meaning
that there is a sequence of vertices $(x_i)_{i = 0,\dots, n}$
with $e_i=(x_{i-1},x_i)$ if $\epsilon_i=1$, and $e_i=(x_i,x_{i-1})$ if
$\epsilon_i=-1$. The inverse of a path $\alpha$ is then defined as
$\alpha^{-1}=e_n^{-\epsilon_n}e_{n-1}^{-\epsilon_{n-1}} \dots
e_1^{-\epsilon_1}$. We also define $begin(\alpha):=x_0$, $end(\alpha):=x_n$.

{\df An Oriented System is a pair $(\Gamma, \sim)$ where $\Gamma$ is
  an oriented graph and $\sim$ is an identification between paths of
  $\Gamma$ such that
\begin{tabular}{cl}
(1) & $\alpha \sim \beta$ implies $begin(\alpha)=begin(\beta)$ and
    $end(\alpha)=end(\beta)$. \\
(2) & $\alpha\alpha^{-1}\sim begin(\alpha)$ for every $\alpha$. \\
(3) & $\alpha\sim\beta$ implies $\alpha^{-1}\sim\beta^{-1}$.\\
(4) & $\alpha\sim\beta$ implies
    $\gamma_1\alpha\gamma_2\sim\gamma_1\beta\gamma_2$\\& for any
    $\gamma_1, \gamma_2$ with $end(\gamma_1)=begin(\alpha)$ and $begin(\gamma_2)=end(\alpha)$.\\
\end{tabular}}

\vspace{0.5cm}

Given a real arrangement of hyperplanes $\AA$, we can associate to it an
oriented system $(\Gamma(\AA),\sim)$. 
Let $\CC(\AA)$ be the vertex set of $\Gamma(\AA)$ and define an
oriented edge pointing from $A$ to $B$ if $A$ is adjacent to $B$
%The chambers of $\AA$ will form the vertex set, and there will be an arrow pointing from $A$ to $B$ if and only if the
%chamber $A$ is adjacent to the chamber $B$ 
(i.e., there is only one hyperplane separating them). 
%Two paths will be identified by $\sim$ 
Identify two paths
if they are both positive (i.e. $\epsilon_i=1 \;\textrm{for all}\; i$), they start at
the same point $v$ and end at the same point $w$, and they both have
minimal length among all positive paths connecting $v$ to $w$.

%The oriented system associated this way with the arrangement $\AA$ will
%be noted $(\Gamma(\AA),\sim)$. 
Note that, geometrically, $\Gamma(\AA)$ is the
1-skeleton of the Salvetti complex.

\subsection{The fundamental group $\pi(\Gamma, \sim )$}

Forgetting orientation of edges, we can view the oriented graph
$\Gamma$ as a 1-complex, and therefore consider its fundamental group
$\pi_1(\Gamma)$ (here we will always deal with connected graphs, so we
do not need to specify a base point). %Although the identification
%$\sim$ was defined on the set of pathes of $\Gamma$, 
The conditions that were required in the
definition of the equivalence relation on paths ensure that $\sim$ induces an equivalence relation on
$\pi_1(\Gamma)$; therefore we can consider the quotient
\begin{displaymath}
\pi(\Gamma, \sim):= \pi_1(\Gamma) / \sim ,
\end{displaymath}
which we will call the \emph{fundamental group} of the oriented system
$(\Gamma, \sim)$.

As a \emph{curiosum}, let us point out that in the case of
$(\Gamma(\AA),\sim)$ the definition of $\sim$ is such that all loops
 $\gamma \in \pi_1(\Gamma(\AA))$ of the form $\gamma=\alpha\beta^{-1}$
for $\alpha, \beta$ minimal positive paths are equivalent to the
trivial loop in $\pi(\Gamma(\AA),\sim)$. These are the loops that by
construction bound some 2-cell in the Salvetti complex. %It will be a
%corollary of some more general later result that in fact 
We will see later (proposition \ref{prop2}) that, in fact,
\begin{displaymath}
\pi(\Gamma(\AA), \sim) \cong \pi_1(\S).
\end{displaymath}

\subsection{Covers of Oriented Systems}

For oriented systems, Paris \cite{Pa1} introduced the following
concept of a cover:

{\df Given two oriented systems $(\Theta, \sim_{\Theta})$ and $(\Psi,
  \sim_{\Psi})$, a morphism of oriented graphs $\rho : \Theta \rightarrow
  \Psi$ is said to be a \textbf{cover} of $(\Psi, \sim_{\Psi})$ if 
\begin{tabular}{cl}
1)& for every vertex $v$ of $\Theta$ and every path $\alpha$ in $\Psi$
   with $begin(\alpha)=\rho(v)$ \\ 
   & there is a unique path $\hat{\alpha}_v$ in $\Theta$ with
   $\rho(\hat{\alpha}_v)=\alpha$ and $begin(\hat{\alpha}_v)=v$. \\
   &This path is called the lift of $\alpha$ at $v$.\\
2)& for any two paths $\alpha,\beta$ in $\Psi$ with
   $begin(\alpha)=begin(\beta)=\rho(v)$, \\
   & if $\alpha\sim\beta$ then $\hat{\alpha}_v\sim\hat{\beta}_v$.
\end{tabular}}

\vspace{0.5cm}
In this setting, Paris also showed a theorem that we recall for later use:

{\thm\label{pa} Consider an oriented system $(\Psi,\sim)$. For each
  subgroup $G$ of $\pi(\Psi,\sim)$ there exists a cover $\rho :
  (\Theta,\sim)\rightarrow(\Psi,\sim)$ with $\rho_\star (\pi(\Theta,\sim))=G$.
}

{\pf See \cite{Pa1}}

\vspace{0.5cm}

In the following we will concentrate on the oriented system
$(\Gamma(\AA),\sim)$ associated with an arrangement, and its
covers. Therefore we will no longer explicitly refer to the arrangement
$\AA$ and simply write $\Gamma$ for $\Gamma(\AA)$.

By lemma 3 of \cite{S}, given two vertices $v_1,v_2$ there is
always a positive minimal path in $\Gamma$ going from $v_1$ to
$v_2$. This motivates the following definition: given a cover
$\rho:(\Theta,\sim)\rightarrow(\Gamma,\sim)$, two vertices $v_1,v_2$
of $\Gamma$ and a vertex $v$ of $\Theta$ with $\rho(v)=v_1$, we will
denote by
\begin{displaymath}
\hat\gamma_v (v_1\rightarrow v_2)
\end{displaymath}
the (unique) positive minimal path of $(\Theta,\sim )$ starting at $v$ and lifting the
(unique) positive minimal path of $(\Gamma,\sim)$ going from $v_1$ to $v_2$.

%%%%%%%%%%%%%%%%%%%%%%%%%%%%%%% DIAGRAM OF SPACES

\section{Diagrams of spaces}

The theory of homotopy colimits of diagrams comes from homological
algebra and  category theory. It was developed by Quillen,
Bousfield, Kan and others (see for example \cite{Q}, \cite{BoKa}), and has now reached remarkable
extension and depth.
In this work we will take the more combinatorial point of view
that Welker, Ziegler and \v Zivaljevi\'c adopted in \cite{WZZ}, where
they developed a useful toolkit for applications of homotopy colimits
in discrete mathematics. We refer to that paper for a comprehensive
and detailed introduction to the subject. 

\subsection{Diagrams of spaces}

A diagram of spaces is a covariant functor
\begin{displaymath}
\mathscr{D}: \II \longrightarrow Top
\end{displaymath}

where $\II$ denotes some small category and $Top$ is the category of
topological spaces and continuous maps.

In our combinatorial setting, $\II$ will always be some poset
$\PP$. Indeed, a poset is a small category if we say that there is an
arrow from $p\in\PP$ to $q\in\PP$ if and only if $p\geq q$ in $\PP$.

In order to simplify notation, we will sometimes write $D_p$ for the
space $\mathscr{D}(p)$, and $f_{p>q}$ for the map $\mathscr{D}(p>q)$.

\subsection{Homotopy colimits}

Given a diagram of spaces, one can consider its homotopy
colimit. Again, we will take here the 'constructive', topological definition given in \cite{WZZ}.

{\df Given a diagram of spaces $\mathscr{D}:\PP\rightarrow Top$, the
  \textbf{homotopy colimit} of $\mathscr{D}$ is defined by
\begin{displaymath}
\HC\mathscr{D} := \bigsqcup_{p\in\PP} D_p \times \Delta(\PP_{\leq p})\bigg/\sim
\end{displaymath}
where the relation $\sim$ is given, for $p>q$, by the following identifications:\\
\begin{displaymath}
\begin{array}{c}
 D_p \times \Delta(\PP_{\leq q}) \hookrightarrow D_p \times
\Delta(\PP_{\leq p})\\ \\
  D_p \times \Delta(\PP_{\leq q}) \stackrel{\scriptsize{(id\times f_{p>q})} }{\longrightarrow}
D_q \times \Delta(\PP_{\leq q})\\
\end{array}
\end{displaymath}
}

\vspace{0.4cm}

As a very small example consider the poset
  with two elements $p$ and $q$ and the order relation given by $p>q$. Let a
  diagram over this poset be given that associates to $p$ the space
  $A$, to q the space $B$, and to the order relation $p>q$ the map
  $f:A\rightarrow B$. Then the homotopy colimit is
\begin{displaymath}
A\times [0,1] \sqcup B\times \{1\} \bigg/ \big\{a\times\{1\} \sim
f(a)\times\{1\} \textrm{ for } a\in A \big\},
\end{displaymath}
hence it is the mapping cylinder of $f$ (see \cite{W} I,5).

\subsection{The simplicial model lemma}

Explicit computations of homotopy colimits can become very long and
cumbersome, even for small examples. Therefore it is convenient to
apply  the tools described in \cite{WZZ} for manipulating a
diagram and bringing it into a more manageable form.

We consider the special case of diagram spaces
$\mathscr{D}(p)$ being homotopy equivalent to order complexes of posets $\QQ_p$. Then it is possible to write the homotopy
colimit as an order complex of an explicitly given poset.

{\df A \textbf{diagram of posets} is a diagram of spaces
  $\mathscr{D}:\PP\rightarrow Top$ such that the spaces assigned to poset elements are order complexes, i.e., there exist posets $\QQ_p$ for $p\in\PP$ such that $\mathscr{D}(p) =
  \Delta(\QQ_p)$ for each $p$.

In this situation, we also can define the \textbf{poset limit} $\PL\mathscr{D}$ of the
diagram of posets $\mathscr{D}$. This is a poset with set of elements
\begin{displaymath}
\PL\mathscr{D} := \bigcup_{p\in\PP} \{p\}\times\QQ_p
\end{displaymath}
and order relations defined by
\begin{displaymath}
(p_1,q_1)\geq(p_2,q_2): \Leftrightarrow 
\left\{\begin{array}{l}
p_1\geq p_2 \quad\textrm{and}\\ 
f_{p_1>p_2}(q_1)\geq q_2 .
\end{array}\right.
\end{displaymath}}

We can now state the \emph{simplicial model lemma:}

{\lm\label{sml} Let $\mathscr{D}$ be a diagram of posets. Then the homotopy colimit of $\mathscr{D}$ is homotopy equivalent to the order complex of the poset limit of $\mathscr{D}$:
\begin{displaymath}
\HC\mathscr{D} \simeq \Delta(\PL\mathscr{D}).
\end{displaymath}}

{\pf  See \cite{BaKo}, note after Corollary 2.11.}

\subsection{The subdivision lemma}

We now present another technical tool that we will need in this
 exposition. It requires some preparatory definitions:

{\df\label{defmin} Given a poset $\mathcal{P}$ and an element
  $P\in\mathcal{P}$, let $\mathcal{P}_{\leq P}$ denote the subposet of $\mathcal{P}$ consisting of
  all elements smaller than or equal to $P$. Moreover,
  $\Sigma^{\vee}(\mathcal{P})$ will denote the
  set of all chains \mbox{$\sigma=P_1<\dots <P_k$} in $\mathcal{P}$, ordered
  by reverse inclusion: $\sigma\geq\sigma '\Leftrightarrow
  \sigma\subseteq\sigma'$.\\
Given a diagram of spaces $\mathcal{D}$ on $\mathcal{P}$, we will
  denote by $\mu^{*}\mathcal{D}$ the inverse image of $\mathcal{D}$
  with respect to the map of posets $\mu:\Sigma^{\vee}(\mathcal{P})\rightarrow
  \mathcal{P}$ that sends each chain to its minimal element.  
The spaces are defined as
 $(\mu^{*}\mathcal{D})(\sigma):=\d(\mu(\sigma))$, and since $\sigma
 >\sigma'$ implies $min(\sigma)\geq min(\sigma')$ in $\FF$, the
  natural choice for the maps
 $(\mu^{*}\d)(\sigma > \sigma')$ is to take the morphisms
 $\d(min(\sigma)\geq min(\sigma'))$.}\\

{\lm\label{lemmin} Let $\mathcal{D}$ be a diagram of spaces. We have
\begin{displaymath}
\HC\mathcal{D}\simeq \HC \mu^{*}\mathcal{D}.
\end{displaymath}}

{\pf Consider an element $P\in\mathcal{P}$. The preimage
  $\mu^{-1}(\mathcal{P_{\leq P}})$ is meet-contractible via
 $\{F\}\in min^{-1}(F)$ (i.e., the meet $\sigma \land \{F\}$ exists for
 all $\sigma\in min^{-1}(\FF_{\leq F})$), hence it is contractible by
 corollary 10.13 of \cite{B}. The claim now follows from lemma 4.7 of \cite{WZZ}. $\square$}

\vspace{0.6cm}
This concludes the preparatory first part of our article. In the
following, we will introduce the diagram models.

\vspace{1cm}  
%%%%%%%%%%%%%%%%%%PART 2%%%%%%%%%%%%%%%%%%%%%%%%%%%%%%%%%%%%

\part{The diagram models.}

Recall that we always have in mind a real arrangement $\AA$ that we
want to study, with its face poset $\FF$ and the associated oriented
system $(\Gamma,\sim)$.

To each cover of the Salvetti complex of $\AA$ we will associate a certain diagram of spaces,
whose homotopy colimit is homotopy equivalent to the cover.
 
%%%%%%%%%%%%%%%%%%%%% SECTION 1%%%%%%%%%%%%%%%%%
\section{The main character}

As all mathematical stories, also ours begins with a definition:

{\df Given a cover of oriented systems $\rho: (\Theta,\sim)\rightarrow
  (\Gamma,\sim)$, $(\Gamma,\sim)$ the oriented system associated to an arrangement $\AA$,  we define a diagram of spaces 
\begin{displaymath}
\mathscr{D}_\rho : \FF \longrightarrow Top
\end{displaymath}
with discrete spaces
\begin{displaymath}
\mathscr{D}_\rho (F) := \{v\in V(\Theta) \vert \rho(v) < F\},
\end{displaymath}
and maps being inclusions
\begin{displaymath}\begin{array}{rccl}
f^\rho_{F_1>F_2}:=\mathscr{D}_\rho(F_1 > F_2):& \mathscr{D}_\rho(F_1)&
\longrightarrow&\mathscr{D}_\rho(F_2)\\
& v & \longmapsto & end\big(\hat{\gamma}_v(\rho(v)\rightarrow\rho(v)_{F_2})\big).  
\end{array}\end{displaymath}}

%The reader is invited to check himself that this object is well-defined.

Let us first make the following observation.

{\prop\label{id} For the identical cover $id:(\Gamma,\sim) \rightarrow (\Gamma,\sim)$ of the oriented system associated with an arrangement $\AA$, the homotopy colimit of $\d$ is homotopy equivalent to the Salvetti complex of $\AA$: 
\begin{displaymath}
\HC\d \simeq \S,
\end{displaymath}
in fact, the order complex of the poset limit $\PL\d$ is identical with the Salvetti complex as simplicial complexes.
}

{\pf
Let us first see, what spaces and maps are in this case:
\begin{displaymath}
\d(F)=\{C\in\CC(\AA)\vert C<F\}
\end{displaymath}
\begin{displaymath}\begin{array}{rccl}
f^{id}_{F_1>F_2}:=\d(F_1 > F_2):& \d(F_1)&
\longrightarrow&\d(F_2)\\
& C & \longmapsto & C_{F_2}\\  
\end{array}\end{displaymath}

The "spaces" are in fact finite sets of points. We can consider them as
order complexes of posets $\QQ_F$ where
$\QQ_F=\{C\in\CC(\AA)\vert\,C<F\}$ as a set, and the order relation is
the empty one: $C_1\geq C_2$ if and only if $C_1 = C_2$.

In this situation we can apply the \emph{simplicial model 
  lemma (\ref{sml})}: $\HC\d$ has the same homotopy type as
  $\Delta(\PL\d)$, which we want to study now. In the following we set
  $W:=\Delta(\PL\d)$.
The vertex set of the simplicial complex $W$ is clearly
  \begin{displaymath}
\mathscr{V}:=\{(F,C)\in\FF\times\CC\vert\,F>C\}.
\end{displaymath} 
For the higher-dimensional simplices we have to look at the chains in
  $\PL\d$. There we have
\begin{displaymath}
(F_1,C_1)\geq(F_2,C_2)\Leftrightarrow\left\{\begin{array}{cl}
1)& F_1\geq F_2 \\
2)& (C_1)_{F_2}\leq C_2\\
\end{array}\right.
\end{displaymath}
Note that in a chain $(F_1,C_1)>(F_2,C_2)>\dots >(F_j,C_j)$ we have
(see \ref{defsal})
$C_i=(C_{i-1})_{F_i}=
%((C_{i-2})_{F_{i-1}})_{F_i}=
\dots=(\dots((C_1)_{F_2})_{F_3}\dots)_{F_i}=(C_1)_{F_2\circ F_3 \dots F_i}=(C_1)_{F_i}$.
A chain in $\PL\d$ is then given by a chain in $\FF$ and a chamber
adjacent to the maximal element of this chain. The rest can be
reconstructed as above. Since all chains are of this form, we can
encode each simplex of $W$ by $\Delta(\phi,C)$, with which we will mean the $d$-simplex given by the $d+1$-chain
$\phi$ and the chamber $C>max(\phi)$. A quick
look at definition \ref{defsal} shows that the evident bijection
$\mathscr{V}\rightarrow V(\S)$ can be linearly extended to all
simplices to get an isomorphism of simplicial complexes.  $\square$}\\

This construction was motivated by a question asked by Michael Falk at
the end of section 2 in \cite{F}, where he presented a model $\hat{S}$
for the complement of complexified line
arrangements. It was constructed starting from the local data of the
arrangements $\mathcal{A}_{X}$ for $X\in\mathcal{L}$ and attaching the
complexes $Sal(\mathcal{A}_X)$ "in the right way". By contraction of
an appropriate subcomplex of $\hat{S}$ he obtained a minimal
2-dimensional complex carrying the homotopy type of the arrangement
complement.\\
The following theorem yields a natural generalization of the
construction of $\hat{S}$ for affine arrangements of any dimension.

{\thm Given a real arrangement $\AA$, define a diagram of spaces
  $\mathcal{E}$ on the poset $\FF^{op}$ by
\begin{displaymath}\begin{array}{c}
\mathcal{E}(F):=Sal(\AA_{\vert F\vert})\\
\mathcal{E}(F_1 >F_2):\textrm{ the natural map of complexes } Sal(\AA_{\vert F_1 \vert}) \hookrightarrow
Sal(\AA_{\vert F_2 \vert}),\\
\textrm{induced by the inclusion }\AA_{\vert
  F_1\vert}\hookrightarrow \AA_{\vert F_2\vert}. 
\end{array}\end{displaymath}
Then we have
\begin{displaymath}
\HC\mathcal{E} \simeq \HC\d \simeq \M.
\end{displaymath}
} 

{\pf In the following, given a diagram $\mathcal{D}$ over a poset
  $\mathcal{P}$ and an element $P\in\mathcal{P}$, let
  $\mathcal{D}[\mathcal{P}_{\leq P}]$ denote the diagram over
  $\mathcal{P}_{\leq P}$ that is naturally induced by $\mathcal{D}$.
  Moreover, $\FF$ will always denote $\FF(\AA)$, the face poset of the
  whole arrangement. 

  Consider a face $F\in\mathcal{F}$. It is easy to see that
  $\mathcal{F}(\AA_{\vert F\vert})=\FF_{\leq F}$. We have
\begin{equation*}
 \mathcal{E}(F)=Sal(\AA_{F})\simeq \HC\d[\FF_{\leq F}].
\end{equation*} 

 We now consider the diagram $\mu^{*}\d$ as defined in
 \ref{defmin}. The subdivision lemma \ref{lemmin} guarantees
\begin{displaymath}
\HC\mu^{*}\d \simeq \HC\d.
\end{displaymath}

 Consider now the map of posets $f:\Sigma^{\vee}(\FF)\rightarrow \FF^{op}$
 defined by $f(\sigma)= max(\sigma)$.
 The homotopy colimit of a diagram over a poset $\mathcal{P}$ with a
 unique minimal element $\hat{0}$ is homotopy equivalent to the space
 attached to $\hat{0}$ (indeed, consider the inclusion
 $\iota\!\!:\!\{\hat{0}\}\!\rightarrow\!\mathcal{P}\;$: since $\hat{0}$ is the
 unique minimal element, the preimage $\iota^{-1}(\mathcal{P}_{\leq
 P})$ is contractible for all $P\in\mathcal{P}$ and we can apply lemma
 4.7 of \cite{WZZ}).
Hence by lemma 4.8 of \cite{WZZ}, there is an homotopy equivalence
\begin{displaymath}
 \HC\mu^{*}\d\simeq \HC f_{\sharp}(\mu^{*}\d),
\end{displaymath}
 where $f_{\sharp}(\mu^{*}\d)$ is defined over the poset
 $\FF^{op}$ by
\begin{equation*}
f_{\sharp}(\mu^{*}\d) (F) :=
\HC(\mu^{*}\d)[f^{-1}(\FF^{op}_{\geq F})],
\end{equation*}
the maps $f_{\sharp}(\mu^{*}\d)(F>F')$ being the natural
inclusions of homotopy colimits.\\

Now we leave it to the reader to check that
 $(\mu^{*}\d)[f^{-1}(\FF^{op}_{\geq F})]= \mu^{*}(\d[\FF_{\leq
 F}])$. For any $F\in\FF^{op}$ we then have an homotopy equivalence by
 lemma \ref{lemmin}:
\begin{equation*}
f_{\sharp}(\mu^{*}\d) (F) \simeq
%%%%%\HC(\mu^{*}\d)[f^{-1}(\FF^{op}_{\geq F})] =
\HC\d[\FF_{\leq F}] \rightarrow \mathcal{E}(F).
\end{equation*}

These maps commute with the morphisms of the diagrams and hence induce
an homotopy equivalence between the homotopy colimits by lemma 4.6 of
\cite{WZZ}.

Putting everything together we have then
\begin{displaymath}
\M\simeq \HC\d\simeq \HC\mu^{*}\d\simeq \HC
f_{\sharp}(\mu^{*}\d) \simeq \HC \mathcal{E}, \qquad
\end{displaymath}
as required. $\square$

}

%%%%%%%%%%%%%%%%%%%%%%%%%%%%%%%COVERING MAPS%%%%%%%%%%%%%%%%%
\section{Covering maps}

We now turn to the general case of a cover
$\rho:(\Theta,\sim)\rightarrow(\Gamma,\sim)$. In fact, all diagrams $\dd$ are diagrams of
posets (see the proof of \ref{id}), and an application of the \emph{simplicial model lemma \ref{sml}} gives 
\begin{displaymath}
\HC\dd \simeq \Delta(\PL\dd)
\end{displaymath}

\vspace{0.3cm}

 For simplicity let us from now on write $\W := \Delta(\PL\dd )$.

The simplicial complex $\W$ has vertex set
\begin{displaymath}
\mathscr{V}_\rho = \{(F,v)\in \FF\times V(\Theta)\vert\, \rho(v)<F\}
\end{displaymath}
and the simplices are chains with respect to the partial order
\begin{displaymath}
(F_1,v_1)\geq(F_2,v_2)\Leftrightarrow\left\{\begin{array}{cl}
1)& F_1 \geq F_2\\
2)& v_2 = end\big(\hat{\gamma}_{v_1}(\rho(v_1)\rightarrow\rho(v_1)_{F_2})\big)\\
\end{array}\right.
\end{displaymath}

A technical lemma has to be worked out here in order to make
the following easier to read.

{\lm For the vertex set $\mathscr{V}_{\rho}$ of $W_{\rho} := \PL\dd$ the following holds:
\begin{displaymath}
\mathscr{V}_\rho = \FF\times V(\Theta) \big/  \sim
\end{displaymath}
where $\sim$ is the smallest equivalence relation generated by
\begin{displaymath}
(F,v)\sim \big(F, end\big(\hat{\gamma}_v(\rho(v)\rightarrow\rho(v)_F) \big)\big).
\end{displaymath}}
{\pf  
Consider the map $j:\V  \rightarrow \FF \times V(\Theta) /  \sim$.

It is surjective, since in the equivalence class of each pair $(F,v)$
there is  some
$(F,\tilde{v})\in \V$, for example the one with
$\tilde{v}=end(\hat{\gamma}_v(\rho(v)\rightarrow \rho(v)_F))$.

The map $j$ would fail to be injective if we could find $(F,v_1),
(F,v_2) \in \V$ with $v_1 \neq v_2$ such that there is  $v\in
V(\Theta)$ with
\begin{displaymath}
v_1 = end\big( \hat{\gamma}_v (\rho(v)\rightarrow\rho(v)_F) \big) = v_2.
\end{displaymath} 
This is impossible by the uniqueness of lifting paths.$\square$}

\vspace{0.7cm}

The map $\rho:\Theta\rightarrow \Gamma$ naturally induces a morphism of
diagrams that we call.
\begin{displaymath}
\lambda:\dd \longrightarrow \d .
\end{displaymath}
On diagram spaces of $\dd$ it is given by a map
\begin{displaymath}\begin{array}{rccl}
\lambda_F :& \dd(F)& \rightarrow & \d(F)\\
           &  v    & \mapsto     &  \rho(v)\\
\end{array}\end{displaymath}
for every $F\in\FF$.

\vspace{0.3cm}
{{\bf Remark:} To check compatibility with the morphisms, consider
  that for each pair  $F_1>F_2$ and $v\in\dd(F_1)$, we have 
\begin{displaymath}\begin{array}{rl}
\lambda_{F_2}(f^\rho_{F_1>F_2}(v))&
=\rho(end(\hat{\gamma}_v(\rho(v)\rightarrow\rho(v)_{F_2})))\\ 
&=\rho(v)_{F_2}=f^{id}_{F_1>F_2}(\rho(v))=f^{id}_{F_1>F_2}(\lambda_{F_1}(v)).\end{array}\end{displaymath}}

The morphism $\lambda$ induces a map between the
homotopy colimits by functoriality, and thus a map 
\begin{displaymath}
\Lambda_\rho : \W \longrightarrow  W.
\end{displaymath}

In fact, $\lambda_\rho$ is a simplicial extension of $\lambda$. The
simplex $\Delta(\phi, v)$ of $\W$ is mapped to $\Delta(\phi,\rho(v))$
in $W$.

Note that the previous considerations can be followed step by step to see
that also a morphism $\eta:(\Theta_1,\sim)\rightarrow (\Theta_2,\sim)$
between two covers $\rho_i :(\Theta_i,\sim)\rightarrow (\Gamma,\sim)$
induces a map $\lambda_\eta : W_{\rho_1}\rightarrow W_{\rho_2}$.

\vspace{0.5cm} 

The goal of this section is to prove
%%%%%%%%%%%%%%%%%%%%%%%%%%%%%%%%PROPOSITION 1%%%%%%%%%%%%
{\prop\label{prop1} The diagram map $\Lambda_\rho: \W \rightarrow W$, induced by a cover of oriented systems $\rho:(\Theta,\sim)\rightarrow (\Gamma,\sim)$, is a topological cover of $W$.}

\vspace{0.3cm}

{\pf The first thing one should ensure when speaking of covering
  spaces is that the base space is connected and locally arcwise
  connected. In our case this follows from the fact that $W$ is finite
  dimensional and locally finite. %If we think of a realization of $W$
  %in some $\R^N$ with simplices being convex hulls of vertices and
  %take a point $P\in W$, there is $\epsilon >0$ such that the open ball
  %$B^N_\epsilon(P)$ does not contain any vertex of $W$ (except $P$ if
  %$P$ happens to be a vertex). Then, denoting by $U$ the star in $W$ of
  %the smallest dimensional simplex that contains $P$, we have a family
  %of neighborhoods of $P$

  Now take $P\in W$ and $\epsilon >0$ such that $B^N_\epsilon (P)$ does
  not contain any vertex of $W$ (except $P$ if $P$ happens to be a
  vertex).  Let $\sigma$ be the smallest dimensional simplex of $W$ containing
  $P$, and let $U$ be the star of $\sigma$. 
%Since $W$ is finite dimensional and locally finite, 
% \begin{displaymath}
% \big\{B^N_{1/n}(P)\cap U \big\}_{n>1/\epsilon}
% \end{displaymath}
%  is a neighborhood base of $P$ consisting of open (in U), arcwise connected sets.
  We have to show that each component of
  $\L^{-1}(B^N_\epsilon (P))$ is mapped homeomorphically to
  $B^N_\epsilon(P)$.
  For this, we will show that $\L^{-1}(U)$ is a disjoint union of
  copies of $U$, each of which is mapped identically to $U$ by $\L$.  
  
\vspace{0.2cm}

  There is a chain $\tilde{\phi}$ and a chamber $\tilde{C}<
  max(\tilde{\phi})$ such that $\sigma=\Delta(\tilde{\phi},\tilde{C})$.
  Defining $\tilde{F}:=max(\tilde{\phi})$, we can write $U$ as
\begin{displaymath}
 U = \bigcup_{\tiny\begin{array}{c}\phi\supseteq\tilde{\phi}\\ C <
 max(\phi)\\ C_{\tilde{F}}=\tilde{C}\end{array}} \Delta(\phi,C),
\end{displaymath} 
  which we can simplify, noting that for each
  $\phi\supseteq\tilde{\phi}$ there is only one $C<max(\phi)$ with
  $C_{\tilde{F}} = \tilde{C}$. This chamber will be called $C(\phi)$. So
  now we can write
\begin{displaymath}
 U = \bigcup_{\tiny\phi\supseteq\tilde{\phi}} \Delta(\phi,C(\phi)),
\end{displaymath}
 and 
\begin{displaymath}
 \L^{-1}(U)=\bigcup_{\tiny\begin{array}{c}\phi\supseteq\tilde{\phi}\\
 v:\rho(v)=C(\phi)\end{array}} \Delta(\phi,v).
\end{displaymath} 

The last notation we now introduce is $v(C,w)$ denoting the vertex of $\Theta$
where a directed path should start, in order to be the lift of the
path $C\rightarrow\rho(W)$ ending in $w$. More precisely, for $w\in
V(\Theta)$ and $C\in V(\Gamma)$ we denote by $v(C,w)$ the (unique)
$v\in V(\Theta)$ such that 
\begin{displaymath}
w=end\big(\hat{\gamma}_v(C\rightarrow \rho(w)) \big).
\end{displaymath}

For a fixed $v'\in \rho^{-1}(\tilde{C})$ define
\begin{displaymath}
W_{v'} := \bigcup_{\phi\supseteq\tilde{\phi}} \Delta(\phi,v(C(\phi),v')).
\end{displaymath}

\begin{description}
%%%%%%%%%%%%%%%%%%%%%%%%%%%%%%CLAIM 1
\item[Claim 1]
\begin{displaymath}
\bigcup_{\tilde{v}\in\rho^{-1}(\tilde{C})} W_{\bar{v}}=\L^{-1}(U)
\end{displaymath}
\item[Proof]Let us first show inclusion from left to right. Take a simplex on the
  left hand side, say $\Delta(\phi, v(C(\phi),\bar{v}))$ for some
  $\phi\supseteq\tilde{\phi}$, $ \bar{v}\in\rho^{-1}(\tilde{C})$. Then
  $\L(\Delta(\phi, v(C(\phi),\bar{v})))=\Delta(\phi,
  \rho(v(C(\phi),\bar{v})))=\Delta(\phi,C(\phi))$, which is a simplex
  of $U$ since by assumption $\phi\supseteq\tilde{\phi}$.\\
  For the other inclusion, take $\Delta(\phi,v)$ for any
  $\phi\supseteq\tilde{\phi}$ and $v\in\rho^{-1}(C(\phi))$. Define
  $\tilde{v}:=end(\hat{\gamma}_v(C(\phi)\rightarrow\tilde{C}))$. By construction, $\rho(\tilde{v})=\tilde{C}$ and
  $v=v(C(\phi),\tilde{v})$: this means that $\Delta(\phi,v)$ is a simplex of
  $W_{\tilde{v}}$.
\vspace{0.5cm}

%%%%%%%%%%%%%%%%%%%%%%%%%%%%CLAIM 2
\item[Claim 2]\emph{Fix $\tilde{v}\in\rho^{-1}(\tilde{C})$, then
  $\L:W_{\tilde{v}}\rightarrow U$ is a homeomorphism.}
\item[Proof]We have already seen that each simplex of $W_{\tilde{v}}$
  maps to a simplex of $U$.\\
  Surjectivity is clear, since given a simplex $\Delta(\phi,C)$ of $U$,
  $\sigma:=\Delta(\phi,v(C,\tilde{v}))$ is a simplex of $W_{\tilde{v}}$ with
  $\L(\sigma)=\Delta(\phi,C)$.\\
  Now it will suffice to show injectivity on
  the vertex sets, since linear extension on each simplex
  will give a PL-map. Take two different
  vertices of $W_{\tilde{v}}$, say $w_1:=(F_1,v(C_1,\tilde{v}))$ and
  $w_2:=(F_2,v(C_2,\tilde{v}))$, and suppose $\L(w_1)=\L(w_2)$.\\
  This means first of all, that $F_1=F_2=:F$ and
  $(C_1)_{F}=(C_2)_{F}=:C$. In order to simplify notation we will write
  $v_1:=v(C_1,\tilde{v})$, $v_2:=v(C_2,\tilde{v})$. Recalling that
  $\V=\FF\times V(\Theta)/\sim$, injectivity is then equivalent to
\begin{displaymath}
end(\hat{\gamma}_{v_1}(C_1\rightarrow C)) =
end(\hat{\gamma}_{v_2}(C_2\rightarrow C)).
\end{displaymath} 
  In order to prove this, we will have to distinguish two cases. By
  assumption, the face $F$ must be by assumption comparable to each element of
  $\tilde{\phi}$, but technical reasons suggest a different treatment
  depending on whether $F$ is bigger or smaller than
  $\tilde{F}=max(\tilde{\phi})$. In any case, one important remark is
  that $C_1$ and $C_2$ are by construction of the form $C(\phi)$ for
  some $\phi$ (say, the $\phi$ giving the smallest simplex of
  $W_{\tilde{v}}$ containing $w_1$ or $w_2$). The definition of
  $C(\phi)$ implies then $(C_1)_{\tilde{F}}=(C_2)_{\tilde{F}}= \tilde{C}$. 
  Before going into further details, it is useful to translate lemma n. 3 of \cite{S} to our language:
\end{description}

{\lm\label{salv} Let $\bar{C}\in V(\Gamma)$, and consider $C\in V(\Gamma)$ such
  that $C<F$ for an $F\in\FF$. Then one positive minimal path from
  $\bar{C}$ to $C$ is the composition of a positive minimal path from
  $\bar{C}$ to $\bar{C}_F$ with a positive minimal path from
  $\bar{C}_F$ to $C$}.

\begin{description}
\item[Case {\boldmath $\tilde{F}>F$}] Since by construction $F>C$, this means
  $\tilde{F}>C$. 
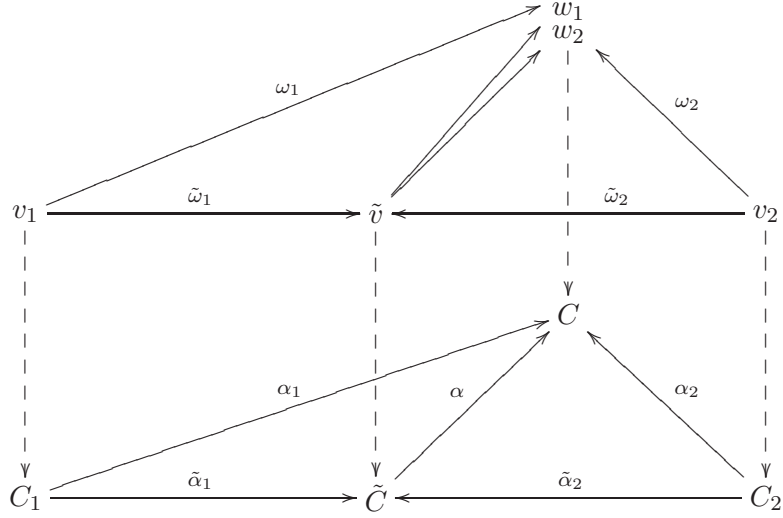
\[
\xymatrix{
                   & &  &  &                      & &{
                   \displaystyle{w_1}\atop \displaystyle{w_2}}\ar@{-->}'[dd][ddd]&&\\
                   & &  &  &                      & &  &  &\\
v_1\ar[rrrr]^{\tilde{\omega}_1}\ar[rrrrrruu]!U^{\omega_1}
   \ar@{-->}[ddd]  & &  &  &\tilde{v}\ar@{-->}[ddd]
                            \ar[rruu]\ar[rruu]!U  & &\ar[ll]&&v_2\ar@{-}[lll]_{\tilde{\omega}_2}
                                                                \ar[lluu]_{\omega_2}
                                                                \ar@{-->}[ddd]\\
                   & &  &  &                      & &C &  &   \\
                   & &  &  &                      & &  &  &   \\
C_1\ar[rrrr]^{\tilde{\alpha}_1}
\ar[rrrrrruu]^{\alpha_1}      
                   & &  &  &\tilde{C}\ar[rruu]^{\alpha}
                                                  & &  &  &C_2\ar[llll]_{\tilde{\alpha}_2}
                                                            \ar[lluu]_{\alpha_2}\\
}
\]

Define 
$\omega_i=\hat{\gamma}_{v_i}(C_i\rightarrow C)$,
%$\omega_2=\hat{\gamma}_{v_2}(C_2\rightarrow C)$, 
$\tilde{\omega}_i=\hat{\gamma}_{v_i}(C_i\rightarrow\tilde{C})$, 
$w_i = end(\omega_i)$,
%$\tilde{\omega}_2=\hat{\gamma}_{v_2}(C_2\rightarrow \tilde{C})$, 
$\tilde{\alpha}_i=\rho(\tilde{\omega}_i)$, for $i=1,2$. 
%$\tilde{\alpha}_2=\rho(\tilde{\omega}_2)$.

Since $\tilde{\alpha}_1$ and $\tilde{\alpha}_2$ are positive minimal paths,
 we obtain a positive minimal path
$\alpha_i:C_i\rightarrow C$ by composing $\tilde{\alpha}_i$ with a
positive minimal path $\alpha :\tilde{C}\rightarrow C$ (see \ref{salv}). So for $i=1,2$:
\begin{displaymath}
 \tilde{\alpha}_i \circ \alpha \sim C_i\rightarrow C
\end{displaymath}
Now we know that $\omega_i =
\hat{\gamma}_{v_i}(\alpha\circ\tilde{\alpha}_i)$, and by construction we
have $end(\tilde{\omega}_1)=\tilde{v}=end(\tilde{\omega}_2)$. With
this we can write for $i=1,2$:
\begin{displaymath}
\qquad end(\omega_i)=end(\hat{\gamma}_{v_i}(\tilde{\alpha}_i \circ \alpha))=end(\hat{\gamma}_{end(\tilde{\omega}_i)}(\alpha))=end(\hat{\gamma}_{\tilde{v}}(\alpha)).
\end{displaymath} 
This clearly does not depend on $i$, and proves the claim in the case $\tilde{F}>F$.

\item[Case {\boldmath $F\leq\tilde{F}$}] In this case both $\tilde{C}$ and $C$ are
  adjacent to $F$. Moreover, recall that, by definition,
  $C=(C_1)_F=(C_2)_F$.
  
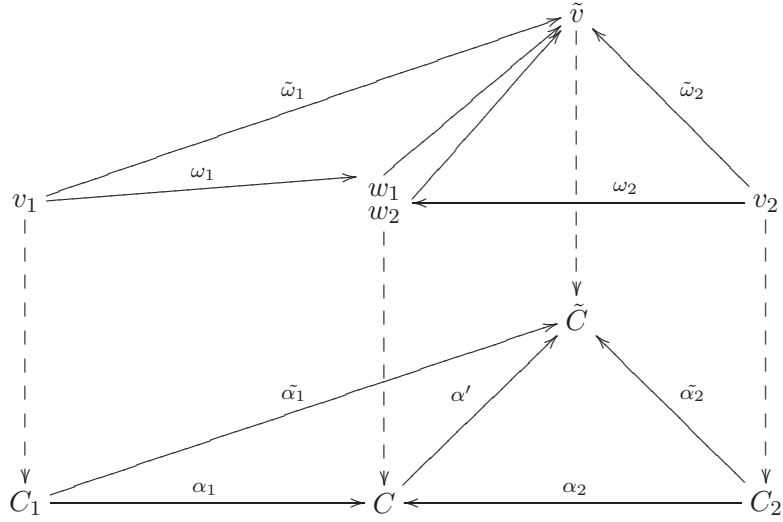
\[
\xymatrix{
                   & &  &  &                      & &\tilde{v}\ar@{-->}'[dd][ddd]
                                                     \ar@{<-}[lldd]!U\ar@{<-}[lldd]!D&&\\
                   & &  &  &                      & &  &  &\\
v_1\ar[rrrrrruu]^{\tilde{\omega}_1}\ar[rrrr]!U^{\omega_1}
   \ar@{-->}[ddd]  & &  &  &{\displaystyle{w_1} \atop \displaystyle{w_2} }\ar@{-->}[ddd]
                                                  & &\ar[ll]  &  &v_2\ar@{-}[lll]_{\omega_2}
                                                                \ar[lluu]_{\tilde{\omega}_2}
                                                                \ar@{-->}[ddd]\\
                   & &  &  &                      & &\tilde{C} &  &   \\
                   & &  &  &                      & &  &  &   \\
C_1\ar[rrrr]^{\alpha_1}
\ar[rrrrrruu]^{\tilde{\alpha_1}}      
                   & &  &  & C\ar[rruu]^{\alpha'}
                                                  & &  &  &C_2\ar[llll]_{\alpha_2}
                                                            \ar[lluu]_{\tilde{\alpha_2}}\\
}
\]

  Let $\omega_i$, $w_i$, $\tilde{\omega}_i$,$\tilde{\alpha}_i$ be defined as
  above and $\alpha'$ be a positive minimal path $C\rightarrow
  \tilde{C}$. Let $\alpha_i=\rho(\omega_i):C_i\rightarrow C$.

  Then by \ref{salv} we know that $\alpha_i\circ\alpha':C_i\rightarrow
  \tilde{C}$ is a positive minimal path, and therefore
  equivalent to $\tilde{\alpha}_i$. This implies
  $\hat{\gamma}_{v_i}(\alpha_i\circ\alpha')\sim\hat{\gamma}_{v_i}(\tilde{\alpha_i})=\tilde{\omega}_i$.
  Consider now $\omega_i^{-1}\circ\tilde{\omega}_i$ for $i=1,2$.\\
  It is easy to see that $\rho(\omega_i^{-1}\circ\tilde{\omega}_i)=
  \alpha^{-1}_i\tilde{\alpha}_i \sim \alpha'$, and because of $end(\omega_i^{-1}\circ\tilde{\omega}_i)=\tilde{v}$,
\begin{displaymath}
\qquad end(\omega_1)=begin(\omega_1^{-1}\circ\tilde{\omega}_1)=v(C,\tilde{v})=begin(\omega_2^{-1}\circ\tilde{\omega}_2)=end(\omega_2).
\end{displaymath} 
  This shows the result in the case $\tilde{F}\geq F$ and concludes
  the proof of \mbox{Claim 2.}
\vspace{0.5cm}
%%%%%%%%%%%%%%%%%%%%%%%%CLAIM 3
\item[Claim 3] For $v_1 \neq v_2 \in \rho^{-1}(\tilde{C})$ we have 
\begin{displaymath}
W_{v_1}\cap W_{v_2} = \emptyset
\end{displaymath}
\item[Proof] Suppose that there is a simplex $w:=\Delta(F,v)\in
  W_{v_1}\cap W_{v_2}$. %w.l.o.g. suppose $\rho(v)<F$. Since $w$ is a
  %vertex, we know $F\in \phi$ for the chain $\phi$ corresponding to
  %the smallest simplex of $W_{v_i}$ containing $w$. 
  The fact that $\sigma$ belongs to both $W_{v_1}$ and $W_{v_2}$
  implies that we can write $v$ in two ways:
\begin{displaymath}
%\begin{array}{rl}
%v&=end\big(\hat{\gamma}_{v(C(\phi),v_1)}(C(\phi)\rightarrow
%C(\phi)_F)\big)\\&\\
%&\qquad=end\big(\hat{\gamma}_{v(C(\phi),v_2)}(C(\phi)\rightarrow C(\phi)_F)\big).\\
%\end{array}
v=v(C(\phi),v_1)=v(C(\phi),v_2)
\end{displaymath}
  This is the same as to say
\begin{displaymath}
\begin{array}{c}
 v_1=end(\hat{\gamma}_v(C(\phi)\rightarrow\rho(v_1)))\\ 
v_2=end(\hat{\gamma}_v(C(\phi)\rightarrow\rho(v_2))),
\end{array}  
 %v(C(\phi),v_1)=v(C(\phi),v)=v(C(\phi),v_2),
\end{displaymath}
  which implies $v_1 = v_2$ since $\rho(v_1)=\rho(v_2)=\tilde{C}$.
  Such a shared vertex $w$ can followly not exist: this concludes the
  proof of Claim 3.
\end{description}

These considerations on the structure of the map $\L$ show that
Proposition \ref{prop1} holds.$\square$
}

%%%%%%%%%%%%%%%%%%%%%%%%%%%%%%%%%%%%%%SECTION 3, FUNDAMENTAL GROUP.

\section{The fundamental group $\pi_1(\HC \DD_\rho)$}

The goal of this section is to obtain more insight into the structure of
the complexes~$\W$. In the end  we will be able to
state and prove Proposition \ref{prop2}, which relates the fundamental
group of a diagram model and of the corresponding oriented system.

The following general observation is at the core of this section:
for any covering $\rho:\Theta\rightarrow\Gamma$, the space $\W=
\Delta(\PL\dd )$ is the \emph{barycentric subdivision of a CW-complex} having
$V(\Theta)$ as vertices, and a $d$-cell $[F,v]$  for each pair
$F\in \FF$, $F\geq\rho(v)$ with $codim(F)=d$.

In fact, if $codim(F)=d$, 
\begin{displaymath}
[F,v]:= \bigcup_{\phi:\,max(\phi)=F} \Delta(\phi,v)
\end{displaymath}
is the barycentric subdivision of a closed $d$-ball.

The cell $[F,v]$ is attached to those vertices $v'$ that can be
written as
\begin{displaymath}
v'=end\big(\hat{\gamma}_v(\rho(v)\rightarrow \rho(v)_{F'})\big)
\end{displaymath} 
with $F'<F$.

\subsection{1-cells}

In particular, between two vertices $v_1$, $v_2$ we have one $1$-cell
for each codimension 1 face $F$ with $F>\rho(v_1)$, $F>\rho(v_2)$,
and for each $v\in V(\Theta)$ such that $F>\rho(v)$ and 
$v_i = end(\hat{\gamma}_v(\rho(v)\rightarrow\rho(v)_{\rho(v_i)}))$ for $i=1,2$.
\[
\xymatrix{
&(F,v)\ar@{-}[dl]\ar@{-}[dr]&\\
(\rho(v_1),v_1)&&(\rho(v_2),v_2)\\
}
\]

As candidates for $v$ we have only $v_1$ and $v_2$: first of all,
since $codim(F)=1$ and $\rho(v)<F$ we only can have have either
$\rho(v)=\rho(v_1)$ or $\rho(v)=\rho(v_2)$. Now, if we suppose
$\rho(v)=\rho(v_1)$, the condition imposes that 
\begin{displaymath}
v_1=end(\hat{\gamma}_v(\rho(v)\rightarrow\rho(v)_{\rho(v_1)}=\rho(v_1)
%[\footnote{See\ref{fenincasina}}]
)),
\end{displaymath} 
which is $v$. Therefore $v_1=v$.\\
For $v_2$ we must then similarly have $v_2 =
end(\hat{\gamma}_{v_1}(\rho(v_1)\rightarrow\rho(v_2)))$, which is true
if and only if there is a directed edge $v_1\rightarrow v_2$ in $\Theta$.

Taking $\rho(v)=\rho(v_2)$ one comes up with a directed edge
$v_2\rightarrow v_1$, and as we can see, there is no third choice.

\vspace{0.5cm}
Summarizing, we can say that we have a $1$-cell between
$[\rho(v_1),v_1]$ and $[\rho(v_2),v_2]$ in~$\W$ for each edge
connecting $v_1$ and $v_2$ in $\Theta$. This amounts to saying that the
$1$-skeleton of $\W$ as a CW-complex is $\Theta$. 

%%%%%%%%%%%%%%%%%%%%%%%%%%%%%%%%%%%%2-CELLS
\subsection{2-cells}

Recall that the fundamental group of a
CW-complex is determined by the 2-skeleton. So we still have to take a
closer look at the 2-cells of $\W$. In the following we will denote
by $\W^i$ the $i$-skeleton of $\W$.

\vspace{0.5cm}
Fix $v\in V(\Theta)$ and $F\in\FF$ with $codim(F)=2$ 
(w.l.o.g. $\rho(v)<F$).
The vertices in the boundary $\partial [F,v]$ are those of the form 
\begin{displaymath}
\big[C, end\big(\hat{\gamma}_v(\rho(v)\rightarrow C)\big)\big],
\end{displaymath}
where $C$ is any chamber adjacent to $F$.

We label the points $end(\hat{\gamma}_v(\rho(v)\rightarrow C))$ as
$v_i$, $i=0,\dots , 2k-1$, assuming w.l.o.g. $v=v_0$, $C_0:=\rho(v_0)$.
So the vertices in $\partial [F,v]$ are now the $[\rho(v_i),
v_i]_{i=0,\dots , 2k-1}$.

\vspace{0.5cm}
Now consider $v_i\neq v_j$, and suppose that in $\partial [F,v]$ an edge
between them exists. This means that $\rho(v_i)$ is adjacent to
$\rho(v_j)$, and there is $F_1$ with $F>F_1 > \rho(v_i)$ and
$F>F_1>\rho(v_j)$ such that this edge can be written as $[F_1,
\tilde{v}]$. To determine whether $\tilde{v}=v_i$ or $v_j$ (which
gives the direction of the edge!) recall that we must have the
following adjacencies in $\PL\dd$:\\

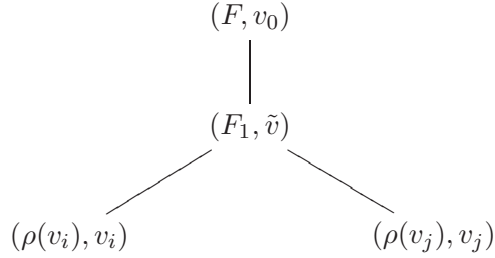
\[
\xymatrix{
&(F,v_0)\ar@{-}[d]&\\
&(F_1,\tilde{v}) \ar@{-}[dl]&\\
(\rho(v_i),v_i) & & (\rho(v_j),v_j)\ar@{-}[ul]\\
}
\]

In particular, $\rho(\tilde{v})=\rho(v_0)_{F_1}=(C_0)_{F_1}$.

In $\partial [F,v]$ we have then one edge for each codimension 1 facet
$F_1$ incident to $F$, and this edge is oriented from the vertex that
projects to the chamber to the same side as $C_0$ w.r.t. $F_1$ to the
other.

\vspace{0.5cm}
We may also view $\partial [F,v]$ as a subgraph of $\W^1$ with an
orientation; namely, a subdivision of $S^1$ with the only one source
($C_0$) and one sink (the chamber opposite to $C_0$, say $C_k$).

Corresponding to it, we have a subgraph in $\Theta$ consisting of two
positive paths from $v_0$ to $v_k$. These are mapped by $\rho$ to two
positive paths $C_0\rightarrow \dots C_i \dots \rightarrow C_k$ with $C_i
<F$ for all $i$. They are both minimal, hence equivalent.

\subsection{The fundamental group}
%\subsection{\boldmath{$\pi_1(\HC \DD_\rho) = \pi(\Theta,\sim)$.}}

We recall the definition of the fundamental group of an oriented
system and have a closer look at the one of $(\Theta,\sim)$, a cover of
$(\Gamma,\sim)$. A loop $\alpha\in\pi(\Theta)$ is equivalent to the
constant loop in $\pi(\Theta,\sim)$ if and only if there are two paths
$\gamma_1$, $\gamma_2$ with $\gamma_1 \sim \gamma_2$ and
$\alpha=\gamma_i\gamma_2^{-1}$. This is precisely the case if
$\rho(\gamma_1)\sim\rho(\gamma_2)$ in $(\Gamma,\sim)$, since $\sim$ on
$\Theta$ is induced through $\rho$ from $(\Gamma,\sim)$ as the
smallest relation such that $\alpha\sim\beta$ if
$begin(\alpha)=begin(\beta)$ and $\rho(\alpha)\sim\rho(\beta)$. We can
state this precisely by saying
\begin{displaymath}
\pi(\Theta,\sim) = \pi(\Theta) \big/ \,\Sigma_\rho,
\end{displaymath}
where $\Sigma_\rho$ is the subgroup generated by the loops of the form
$\gamma_1\gamma_2^{-1}$ with $\rho(\gamma_1)\sim\rho(\gamma_2)$.

As a final remark, note that $\rho(\gamma_1)\sim\rho(\gamma_2)$ in
$(\Gamma,\sim)$ if and only if the corresponding paths are homotopic in $\S$.

\vspace{0.5cm}
Now consider $\pi_1(\W)$, the fundamental group of $\HC\dd$. By
standard tools of algebraic topology (see for ex. \cite{M} ch. 7
Thm. 2.1) we can define
\begin{displaymath}\begin{array}{rl}
\sigma_\rho:= & \textrm{\emph{the smallest normal subgroup of
}}\pi(\W)\textrm{\emph{ containing}}\\
&\{\partial [F,v]\vert\, codim(F)=2, F>\rho(v)\},
\end{array}\end{displaymath}
and say
\begin{displaymath}
\pi_1(\W)=\pi(\W^1)\big/ \, \sigma_\rho.
\end{displaymath}

\vspace{0.5cm}
Now consider the natural inclusion $i:\Theta\rightarrow\W$.
We are now ready to prove

{\lm\label{mainl}
The induced map
\begin{displaymath}
i_\star : \Sigma_\rho \rightarrow \sigma_\rho
\end{displaymath}
\hspace{9cm}is an isomorphism.
}

\vspace{0.3cm}

{\pf First note that $\Sigma_\rho$ is a normal subgroup of
  $\pi(\Theta)$:  for each $\alpha\in\pi(\Theta)$ and
  $\gamma_1\gamma_2^{-1}\in\Sigma_\rho$, we have
  $\alpha\gamma_1\gamma_2^{-1}\alpha^{-1}=
  \alpha\gamma_1\alpha^{-1}\alpha\gamma_2^{-1}\alpha^{-1}\in
  \Sigma_\rho$ since 
  $\rho(\alpha\gamma_1\alpha^{-1})=\rho(\alpha)\rho(\gamma_1)\rho(\alpha)^{-1}\sim\rho(\gamma_1)\sim\rho(\gamma_2)=\rho(\alpha)^{-1}\rho(\gamma_2)\rho(\alpha)=
  \rho(\alpha^{-1}\gamma_2^{-1}\alpha)$.

Therefore, $\Sigma_\rho$ is the smallest normal subgroup containing
$\gamma_1\gamma_2^{-1}$ for each pair $(\gamma_1,\gamma_2)$ with
$\rho(\gamma_1)\sim\rho(\gamma_2)$.

We have already seen that each $i_\star^{-1}(\partial [F,v])$ is always of
the form $\gamma_1\gamma_2^{-1}$ with
$\rho(\gamma_1)\sim\rho(\gamma_2)$.
Let $\omega_i:=\rho(\gamma_i)$: by definition of $\omega_1\sim\omega_2$, they are constructed by concatenation and
conjugation from positive, minimal paths. Since $\Sigma_\rho$ is a
normal subgroup, it remains to show that \emph{for each pair
$(\gamma_1,\gamma_2)$ with $\omega_1\sim\omega_2$ and  $\omega_i$ positive and minimal, there are
$\tilde{\omega}_1$, $\tilde{\omega}_2$, $\alpha$, $F$, $v$ with
$\rho(v)=end(\alpha)$ and such that
\begin{displaymath}
\omega_i\sim\alpha\tilde{\omega}_i\textrm{ for } i=1,2
\end{displaymath} 
and 
\begin{displaymath}
i\big(\hat{\gamma}_v(\tilde{\omega}_1)\hat{\gamma}_v(\tilde{\omega}_2)^{-1}\big)=\partial [F,v].
\end{displaymath}}
First notice that given $\omega_1\sim\omega_2$ positive and minimal in
$\Gamma$ and $F>end(\omega_1)=end(\omega_2)$ with $codim(F)=2$, by
\ref{salv} there is a chamber $C>F$, positive and minimal paths
$\tilde{\omega}_1$, $\tilde{\omega}_2$, $\alpha$ such that
$\omega_i\sim\alpha\tilde{\omega_i}$ for $i=1,2$ and
$end(\alpha)=C$. Moreover, each vertex of $\tilde{\omega}_i$ is
adjacent to $F$, but no vertex of $\alpha$ other than $C$ is adjacent
to $F$.

So given a pair $\gamma_1\sim\gamma_2$ in $\Theta$ with $\omega_1$,
$\omega_2$ positive and minimal as above, let
$v':=begin(\gamma_1)=begin(\gamma_2)$. Applying the previous
considerations we can write
$\rho(\gamma_1\gamma_2^{-1})=\omega_1\omega_2^{-1}=\alpha\tilde{\omega}_1\alpha^{-1}\tilde{\omega}_2^{-1}\sim\tilde{\omega}_1\tilde{\omega}_2^{-1}$.\\
Consider then $\tilde{\gamma}_i:=\hat{\gamma}_v(\tilde{\omega}_1)$, where
$v:=end(\hat{\gamma}_{v'}(\alpha))$.
We have shown that
$\gamma_1\gamma_2^{-1}\sim\tilde{\gamma}_1\tilde{\gamma}_2^{-1}$:  it
is now straightforward to check $i_\star (\tilde{\gamma}_1\tilde{\gamma}_2^{-1})=\partial [F,v]$.

This concludes the proof of the lemma. $\square$
}

\vspace{0.5cm}
The following proposition is now immediate:

{\prop\label{prop2}
For a cover of oriented systems
$\rho:(\Theta,\sim)\rightarrow(\Theta,\sim)$, the fundamental group of
the homotopy colimit of the diagram model $\dd$ is isomorphic to the
fundamental group of the oriented system $(\Theta,\sim)$:
\begin{center}$\pi_1(\HC \DD_\rho) \cong \pi(\Theta,\sim).$\end{center}
}

\section{Results and examples}

%We are now able to formulate and prove some results that follow from
%the previous study of the diegram models.

\subsection{The main result}
%{\thm For each cover $\mu : \tilde{\MM}\rightarrow\M$ there is a
%  diagram model $\dd$ such that $\HC\dd$ and $\tilde{\MM}$ are isomorphic
%  covers of $\M$. 
%}

%\vspace{0.4cm}
%{\pf
%By \ref{prop2} we have an isomorphism $\varphi :
%\pi(\Gamma,\sim)\rightarrow\pi_1(\HC\d)$. Since $\HC\d = \S
%\simeq \M$ by \ref{id} and the homotopy equivalence
%$\Phi:\S\rightarrow \M$ given in \cite{S}, we can consider the preimage
%$G:=\varphi^{-1}(\mu_\star(\pi_1(\tilde{\MM})))$ of the fundamental group
%of $\tilde{\MM}$ in $\pi(\Gamma,\sim)$.

%With theorem \ref{pa} we know that there is a cover
%$\rho:(\Theta,\sim)\rightarrow (\Gamma,\sim)$ with
%$\rho_\star(\pi(\Theta,\sim))= G$. Again by \ref{prop2} we have an
%isomorphism $\varphi_\rho : \pi(\theta,\sim) \rightarrow
%\pi_1(\HC\dd)$. These isomorphisms come naturally from the inclusion of the
%oriented systems as 1-skeleton of the CW-version of the homotopy
%colimits. Therefore the following diagrams commute
%\[
%\xymatrix{(\Theta,\sim)\ar[r]^\iota\ar[d]^\rho&\W\ar[d]^{\L}& 
%          \pi(\Theta,\sim)\ar[r]^{\varphi_\rho}\ar[d]^{\rho_\star}&\pi_1(\W)\ar[d]^{(\L)_\star}\\
%          (\Gamma,\sim)\ar[r]^\iota&W&
%          \pi(\Gamma,\sim)\ar[r]^\varphi& \pi_1(W)
%}
%\]

%and
%$\mu_\star(\pi_1(\tilde{\MM}))=\Phi_\star\varphi\rho_\star(\pi(\Theta,\sim))=(\Phi\L)_\star\varphi_\rho
%(\pi_1(\W))\cong (\Phi\L)_\star(\pi_1(\W))$.
%Therefore the cover $\Phi\L : \W \rightarrow \M$ is isomorphic to
%$\mu: \tilde{M}\rightarrow \M$.
%}

Combining the results of pur studies we are now able to formulate and
prove our main theorem.

{\thm For any topological cover $r : S'\rightarrow\S$ of the Salvetti
  complex of an arrangement $\AA$, there exists a cover of oriented
  systems $\rho:(\Theta,\sim)\rightarrow(\Gamma,\sim)$ such that the
  homotopy colimit of the associated diagram of spaces $\dd$ is
  homotopy equivalen to $S'$. Moreover, the poset limit of $\dd$ is
  isomorphic to $S'$ as a covering space of $\S$.
}

\vspace{0.4cm}
{\pf
By \ref{prop2} we have an isomorphism $\varphi :
\pi(\Gamma,\sim)\rightarrow\pi_1(\HC\d)$. Since $\HC\d \simeq \S $ by \ref{id},  we can consider the preimage
$G:=\varphi^{-1}(r_\star(\pi_1(S')))$ of the fundamental group
of $S'$ in $\pi(\Gamma,\sim)$.

Theorem \ref{pa} gives a cover
$\rho:(\Theta,\sim)\rightarrow (\Gamma,\sim)$ with
$\rho_\star(\pi(\Theta,\sim))= G$. Again by \ref{prop2} we have an
isomorphism $\varphi_\rho : \pi(\Theta,\sim) \rightarrow
\pi_1(\HC\dd)$. These isomorphisms come naturally from the
inclusion $\iota$ of the
oriented systems as 1-skeleton of the CW-version of the homotopy
colimits. Therefore the following diagrams commute
\[
\xymatrix{(\Theta,\sim)\ar[r]^\iota\ar[d]^\rho&\W\ar[d]^{\L}& 
          \pi(\Theta,\sim)\ar[r]^{\varphi_\rho}\ar[d]^{\rho_\star}&\pi_1(\W)\ar[d]^{(\L)_\star}\\
          (\Gamma,\sim)\ar[r]^\iota&W&
          \pi(\Gamma,\sim)\ar[r]^\varphi& \pi_1(W)
}
\]

and
$r_\star(\pi_1(S'))=\varphi\rho_\star(\pi(\Theta,\sim))=(\L)_\star\varphi_\rho
(\pi_1(\W))\cong (\L)_\star(\pi_1(\W))$.
Therefore the cover $\L : \W \rightarrow \S$ is isomorphic to
$r: S'\rightarrow \S$.$\square$
}

\vspace{0.5cm}
An immediate consequence of the theorem is

{\crl Let $\hat{\rho}:(\hat{\Gamma},\sim)\rightarrow (\Gamma,\sim)$
  denote the universal cover of $(\Gamma,\sim)$. Then $W_{\hat{\rho}}=\Delta (\PL\mathcal{D}_{\hat{\rho}} )$ is the universal
  cover of $\S$. }

{\pf We prove universality. Take any cover $r:S'\rightarrow\S$;
  we have to show that there is a morphism of covers
  $m:W_{\hat{\rho}}\rightarrow S'$.
  By the theorem, we know that there is a cover
  $\rho:(\Theta,\sim)\rightarrow (\Gamma,\sim)$ with $\W \cong
  S'$ as a cover. Universality of $(\hat{\Gamma},\sim)$ implies
  the existence of a morphism $\mu :
  (\hat{\Gamma},\sim)\rightarrow(\Theta,\sim)$, and this induces a
  morphism of diagrams
  $\lambda_\mu:\DD_{\hat{\rho}}\rightarrow\dd$. By functoriality,
  we have $\Lambda_\mu:W_{\hat{\rho}}\rightarrow\W$, which gives the
  required morphism, as in the following diagram. 
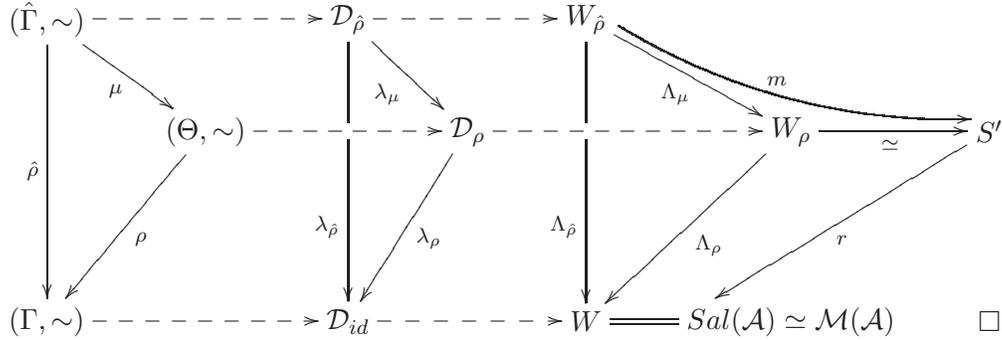
\[ 
\xymatrix{
(\hat{\Gamma},\sim)\ar[dr]_{\mu}\ar[ddd]_{\hat{\rho}}\ar@{-->}[rr]
 &&\DD_{\hat{\rho}}\ar[dr]_{\lambda_\mu}\ar'[d][ddd]_{\lambda_{\hat{\rho}}}\ar@{-->}[rr]&&
W_{\hat{\rho}}\ar[dr]_{\Lambda_\mu}\ar'[d][ddd]_{\Lambda_{\hat{\rho}}}\ar@/_1pc/@<1ex>[drr]!L^m&&\\
%&&&&&&\\
&(\Theta,\sim)\ar[ddl]^{\rho}\ar@{-->}[rr]&&\dd\ar[ddl]^{\lambda_\rho}\ar@{-->}[rr]&&\W\ar[ddl]^{\L}\ar[r]_{\simeq}&S'\ar[ddl]!L^{r}\\
&&&&&&\\
(\Gamma,\sim)\ar@{-->}[rr]&&\d\ar@{-->}[rr]&& W\ar@{=}[r]&\S\simeq\M& \square
}\]
}

\vspace{0.5cm}

As an illustration of the result, we want to work out an easy but
instructive example. 

\subsection{An example}

Consider the arrangement given by one point $P\in\R$. The space $\R$ is divided by $P$
 in two chambers $A$ and $B$. It is easy to write
down the face poset $\FF$ and the oriented system $(\Gamma,\sim)$ as follows:

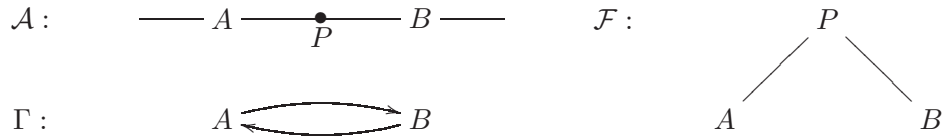
\[
\xymatrix{
\AA:&\ar@{-}[r]&A\ar@{-}[rr]_{\displaystyle P}&\bullet&B\ar@{-}[r]&&\FF:& &P\ar@{-}[dl]\ar@{-}[dr]&\\ 
\Gamma: &&A \ar@/^/[rr]&& B\ar@/^/[ll]&&    &A& &B \\}
\]\\

The complexification of $\AA$ is the arrangement given by a point in
the complex plane. The complement $\M$ is then homotopy equivalent to $S^1$,
and its universal cover hence is $\R$. We will now see how the
diagram models come to this conclusions.

First consider $\d$, a diagram on the poset $\FF$. The 'space'
associated to $A$ and $B$ is one point: $\d(A)=\{A\}$,
$\d(B)=\{B\}$. To the element $P$ we attach a space
consisting of two points, one for each chamber adjacent to $P$. To
keep track of which point is associated to which chamber, we will write
$\d(P)=\{A,B\}$.
The maps between spaces are in this case trivial, but let us explain
where they come from:\\

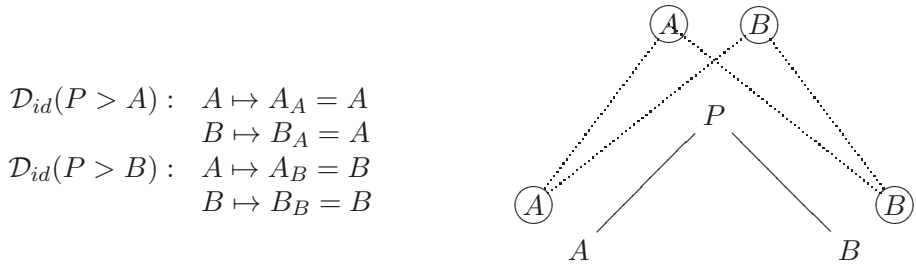
\begin{tabular}{cp{1cm}c}
$\begin{array}{rl}
\d(P>A):& A\mapsto A_A = A\\
        & B\mapsto B_A =A \\
\d(P>B):& A\mapsto A_B =B \\
        & B\mapsto B_B =B \\
\end{array}$%\end{displaymath}
 &&  \xy <0.6cm,0cm>:
(-2,-2)*+{A}; (1,1)*+{P};**\dir{-},(4,-2)*+{B};**\dir{-}; (0,3)*+{A}*\cir{};(-3,-1)*+{A}*\cir{};**\dir{..},
(2,3)*+{B}*\cir{};**\dir{..},
(5,-1)*+{B}*\cir{};**\dir{..}, (0,3);**\dir{..}
\endxy  \\
\end{tabular}
\\
\vspace{0.5cm}

We now have to look at the
universal cover of $(\Gamma,\sim)$. Since $\S$ has no 2-cells, the
identification $\sim$ on $\Gamma$ (hence on $\hat{\Gamma}$, too) is
empty.
Let us then write $(\hat{\Gamma},\sim)$ as follows: 
\[\xymatrix{
%\Gamma: &A \ar@/^/[rr]&& B\ar@/^/[ll]&\FF:& &P\ar@{-}[dl]\ar@{-}[dr]&\\
%&&&                                  &    &A& &B \\
\hat{\Gamma}:&\ar@{-->}[r]&A_{-1}\ar[r]&B_{-1}\ar[r]&A_0\ar[r]&B_0\ar[r]&A_1\ar@{-->}[r]&\\
}\]

where $\r: (\hat{\Gamma},\sim)\rightarrow (\Gamma,\sim)$ is clearly
defined by $\r(A_i)=A$, $\r(B_i)=B$ for all $i$.\\
Writing down the diagram $\DD_{\r}$ we have to keep in mind that the
space associated to an element $F\in\FF$ has as many disconnected
points as there are vertices in $\hat{\Gamma}$ that project to a
chamber adjacent to $F$. So we have $\DD_{\r}(A)=\{A_i \vert\,
i\in\mathbb{Z}\}$, $\DD_{\r}(B)=\{B_i \vert\,
i\in\mathbb{Z}\}$, $\DD_{\r}(P)=\DD_{\r}(A)\cup\DD_{\r}(B)$.
For the maps one has to take care of how paths are lifted. Let us work
out some special case and write down the diagram in the same fashion as above:

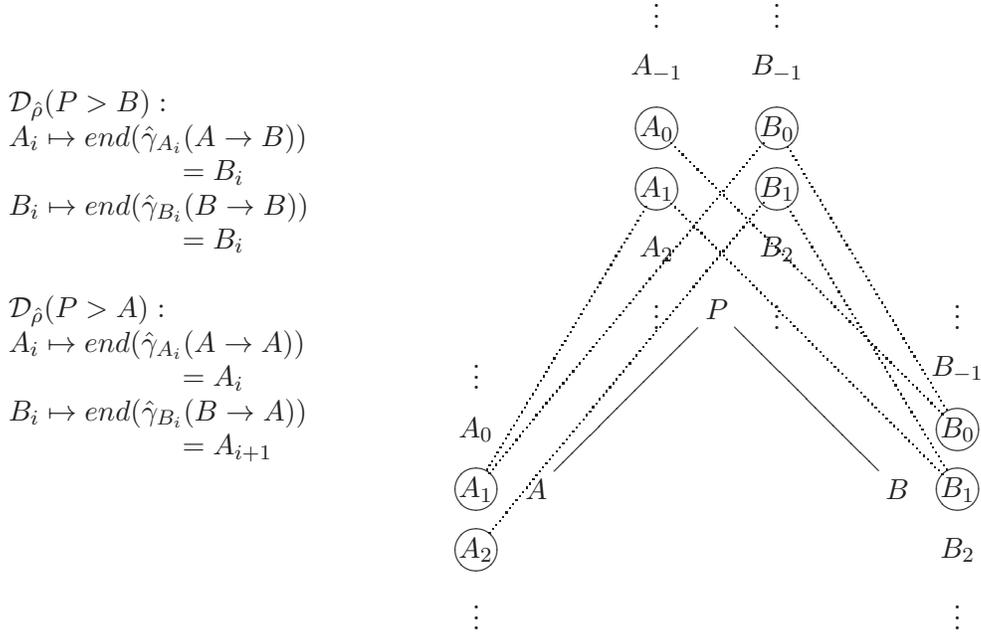
\begin{tabular}{cp{1cm}c}
$\begin{array}{l}
\DD_{\r}(P>B): \\ A_i \mapsto end(\hat{\gamma}_{A_i}(A \rightarrow B))\\ 
\qquad\qquad\qquad =B_i\\
B_i \mapsto end(\hat{\gamma}_{B_i}(B \rightarrow B))\\
\qquad\qquad\qquad =B_i\\

\qquad\\
\DD_{\r}(P>A): \\ A_i \mapsto end(\hat{\gamma}_{A_i}(A \rightarrow A))\\
\qquad\qquad\qquad =A_i\\

B_i \mapsto end(\hat{\gamma}_{B_i}(B \rightarrow A))\\
\qquad\qquad\qquad =A_{i+1}\\

\qquad\\
\qquad\\
\quad\\
\quad\\
\quad\\
%\quad\\
%\quad\\
%\quad\\
%\quad\\
\end{array}$
&&
\xy<0.8cm,0cm>:
(-2,-2)*+{A}; (1,1)*+{P};**\dir{-},(4,-2)*+{B};**\dir{-};
(-3,-3)*+{A_2}*\cir{};(2,3)*+{B_1}*\cir{};**\dir{..}, 
(5,-2)*+{B_1}*\cir{};**\dir{..}, 
(0,3)*+{A_1}*\cir{};**\dir{..}, 
(-3,-2)*+{A_1}*\cir{};**\dir{..},
(2,4)*+{B_0}*\cir{};**\dir{..},
(5,-1)*+{B_0}*\cir{};**\dir{..},
(0,4)*+{A_0}*\cir{};**\dir{..},
(0,6)*+{\vdots},(0,5)*+{A_{-1}},(0,2)*+{A_2},(0,1)*+{\vdots},
(5,0)*+{B_{-1}},(5,1)*+{\vdots},(5,-3)*+{B_2},(5,-4)*+{\vdots},
(2,6)*+{\vdots},(2,5)*+{B_{-1}},(2,2)*+{B_2},(2,1)*+{\vdots},
(-3,-4)*+{\vdots},(-3,0)*+{\vdots},(-3,-1)*+{A_0},
\endxy
\end{tabular}\\
\vspace{0.5cm}

By lemma \ref{sml}, we now only have to write down the poset
$\PL\DD_{\r}$. The order relation is such that for $F_i\in\FF$
and $v_i \in V(\hat{\Gamma})$ we have 
$(F_1,v_1)\geq(F_2,v_2)$ if and only if $F_1 \geq F_2$
and $ v_2 = end(\hat{\gamma}_{v_1}(\rho(v_1)\rightarrow\rho(v_1)_{F_2}))$.

In our case, this means that the dotted lines in the above picture are
yet a piece of the Hasse diagram of $\PL\DD_{\r}$, which we can
redraw in a more readable way as

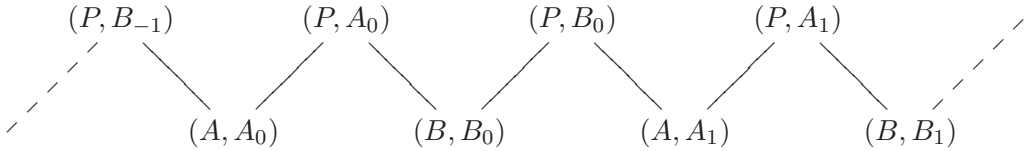
\[\xy<1.5cm,0cm>:
(-3,-1)*+{};(-2,0)*+{(P,B_{-1})};**\dir{--},
%(-3,-1)*+{(B,B_{-1})};**\dir{-},(-2,1)*+{(P,B_{-1})};**\dir{-},
(-1,-1)*+{(A,A_0)};**\dir{-},(0,0)*+{(P,A_0)};**\dir{-},
(1,-1)*+{(B,B_0)};**\dir{-},(2,0)*+{(P,B_0)};**\dir{-},
(3,-1)*+{(A,A_1)};**\dir{-},(4,0)*+{(P,A_1)};**\dir{-},
(5,-1)*+{(B,B_1)};**\dir{-},(6,0)*+{};**\dir{--},
%(7,-1)*+{(A,A_{2})};**\dir{-},(8,0)*+{};**\dir{--},
\endxy\]
%&(P,A_{-1})\ar@{-}[dl]\ar@{-}[dr]&       &(P,B_{-1})\ar@{-}[dl]\ar@{-}[dr]&
%&(P,A_{0})\ar@{-}[dl]\ar@{-}[dr]&        &(P,B_{0})\ar@{-}[dl]\ar@{-}[dr]&
%&(P,A_{-1})\ar@{-}[dl]\ar@{-}[dr]&       &(P,A_{-1})\ar@{-}[dl]\ar@{-}[dr]&\\
%&                                &(B,B_{-1})
%&                                &(A,A_0)   
%&&(B,B_0)
%&&(A,A_1)
%&&(B,B_1)
%&&(A,A_2)\\
\nopagebreak[4]
It is now clear that $\HC\DD_{\rho}\cong\Delta(\PL\DD_{\r})\simeq
\R$, as required.

\end{document}